# NONLINEAR ESTIMATION FOR LINEAR INVERSE PROBLEMS WITH ERROR IN THE OPERATOR[1]

By Marc Hoffmann and Markus Reiss

*University of Marne-la-Vallée and University of Heidelberg*

We study two nonlinear methods for statistical linear inverse problems when the operator is not known. The two constructions combine Galerkin regularization and wavelet thresholding. Their performances depend on the underlying structure of the operator, quantified by an index of sparsity. We prove their rate-optimality and adaptivity properties over Besov classes.

## 1. Introduction.

*Linear inverse problems with error in the operator.* We want to recover $f \in L^2(\mathcal{D})$, where $\mathcal{D}$ is a domain in $\mathbb{R}^d$, from data

$$g_\varepsilon = Kf + \varepsilon \dot{W}, \quad (1.1)$$

where $K$ is an unknown linear operator $K : L^2(\mathcal{D}) \to L^2(\mathcal{Q})$, $\mathcal{Q}$ is a domain in $\mathbb{R}^q$, $\dot{W}$ is Gaussian white noise and $\varepsilon > 0$ is the noise level. We do not know $K$ exactly, but we have access to

$$K_\delta = K + \delta \dot{B}. \quad (1.2)$$

The process $K_\delta$ is a blurred version of $K$, polluted by a Gaussian operator white noise $\dot{B}$ with a noise level $\delta > 0$. The operator $K$ acting on $f$ is unknown and treated as a nuisance parameter. However, preliminary statistical inference about $K$ is possible, with an accuracy governed by $\delta$. Another equivalent approach is to consider that for experimental reasons we never have access to $K$ in practice, but rather to $K_\delta$. The error level $\delta$ can be linked to the accuracy of supplementary experiments; see Efromovich and Koltchinskii [11] and the examples below. In most interesting cases $K^{-1}$

Received December 2004; revised April 2007.
[1]Supported in part by the European research training network Dynstoch.

is not continuous and the estimation problem (1.1) is ill-posed (e.g., see Nussbaum and Pereverzev [16] and Engl, Hanke and Neubauer [12]).

The statistical model is thus given by the observation $(g_\varepsilon, K_\delta)$. Asymptotics are taken as $\delta, \varepsilon \to 0$ simultaneously. In probabilistic terms, observable quantities take the form

$$\langle g_\varepsilon, k \rangle := \langle Kf, k \rangle_{L^2(\mathcal{Q})} + \varepsilon \langle \dot{W}, k \rangle \qquad \forall k \in L^2(\mathcal{Q})$$

and

$$\langle K_\delta h, k \rangle := \langle Kh, k \rangle_{L^2(\mathcal{Q})} + \delta \langle \dot{B} h, k \rangle \qquad \forall (h, k) \in L^2(\mathcal{D}) \times L^2(\mathcal{Q}).$$

The mapping $k \in L^2(\mathcal{Q}) \mapsto \langle \dot{W}, k \rangle$ defines a centered Gaussian linear form, with covariance

$$\mathbb{E}[\langle \dot{W}, k_1 \rangle \langle \dot{W}, k_2 \rangle] = \langle k_1, k_2 \rangle_{L^2(\mathcal{Q})}, \qquad k_1, k_2 \in L^2(\mathcal{Q}).$$

Likewise, $(h, k) \in L^2(\mathcal{D}) \times L^2(\mathcal{Q}) \mapsto \langle \dot{B} h, k \rangle$ defines a centered Gaussian bilinear form with covariance

$$\mathbb{E}[\langle \dot{B} h_1, k_1 \rangle \langle \dot{B} h_2, k_2 \rangle] = \langle h_1, h_2 \rangle_{L^2(\mathcal{D})} \langle k_1, k_2 \rangle_{L^2(\mathcal{Q})}.$$

If $(h_i)_{i \geq 1}$ and $(k_i)_{i \geq 1}$ form orthonormal bases of $L^2(\mathcal{D})$ and $L^2(\mathcal{Q})$, respectively—in particular, we will consider hereafter wavelet bases, the infinite vector $(\langle \dot{W}, k_j \rangle)_{j \geq 1}$ and the infinite matrix $(\langle \dot{B} h_i, k_j \rangle)_{i,j \geq 1}$ have i.i.d. standard Gaussian entries. Another description of the operator white noise is given by stochastic integration using a Brownian sheet, which can be interpreted as a white noise model for kernel observations; see Section 2 below.

*Main results.* The interplay between $\delta$ and $\varepsilon$ is crucial: if $\delta \ll \varepsilon$, one expects to recover model (1.1) with a known $K$. On the other hand, we will exhibit a different picture if $\varepsilon \ll \delta$. Even when the error $\varepsilon$ in the signal $g_\varepsilon$ dominates $\delta$, the assumption $\delta \neq 0$ has to be handled carefully. We restrict our attention to the case $\mathcal{Q} = \mathcal{D}$ and nonnegative operators $K$ on $L^2(\mathcal{D})$.

We first consider a linear estimator based on the Galerkin projection method. For functions in the $L^2$-Sobolev space $H^s$ and suitable approximation spaces, the linear estimator converges with the minimax rate $\max\{\delta, \varepsilon\}^{2s/(2s+2t+d)}$, where $t > 0$ is the degree of ill-posedness of $K$.

For spatially inhomogeneous functions, like smooth functions with jumps, linear estimators cannot attain optimal rates of convergence; see, for example, Donoho and Johnstone [10]. Therefore we propose two nonlinear methods by separating the two steps of Galerkin *inversion* $(I)$ and adaptive *smoothing* $(S)$, which provides two strategies:

Nonlinear Estimation I: $\qquad (g_\varepsilon, K_\delta) \xrightarrow{(I)} \hat{f}_{\delta,\varepsilon}^{\mathrm{lin}} \xrightarrow{(S)} \hat{f}_{\delta,\varepsilon}^I,$

Nonlinear Estimation II: $\qquad (g_\varepsilon, K_\delta) \xrightarrow{(S)} (\hat{g}_\varepsilon, \hat{K}_\delta) \xrightarrow{(I)} \hat{f}_{\delta,\varepsilon}^{II},$

where $\hat{f}_{\delta,\varepsilon}^{\mathrm{lin}}$ is a preliminary and undersmoothed linear estimator. We use a Galerkin scheme on a high-dimensional space as inversion procedure ($I$) and wavelet thresholding as adaptive smoothing technique ($S$), with a level-dependent thresholding rule in Nonlinear Estimation I and a noise reduction in the operator by entrywise thresholding of the wavelet matrix representation in Nonlinear Estimation II. To our knowledge, thresholding for the operator is new in a statistical framework.

From both mathematical and numerical perspectives, the inversion step is critical: we cannot choose an arbitrarily large approximation space for the inversion, even in Nonlinear Estimation II. Nevertheless, both methods are provably rate-optimal (up to a log factor in some cases for the second method) over a wide range of (sparse) nonparametric classes, expressed in terms of Besov spaces $B_{p,p}^s$ with $p \leq 2$.

*Organization of the paper.* Section 2 discusses related approaches. The theory of linear and nonlinear estimation is presented in Sections 3 to 5. Section 6 discusses the numerical implementation. The proofs of the main theorems are deferred to Section 7 and the Appendix provides technical results and some tools from approximation theory.

## 2. Related approaches with error in the operator.

*Perturbed singular values.* Adhering to a singular-value decomposition approach, Cavalier and Hengartner [3] assume that the singular functions of $K$ are known, but not its singular values. Examples include convolution operators. By an oracle-inequality approach, they show how to reconstruct $f$ efficiently when $\delta \leq \varepsilon$.

*Physical devices.* We are given an integral equation $Kf = g$ on a closed boundary surface $\Gamma$, where the boundary integral operator

$$Kh(x) = \int_\Gamma k(x,y)h(y)\sigma_\Gamma(dy)$$

is of order $t > 0$, that is, $K : H^{-t/2}(\Gamma) \to H^{t/2}(\Gamma)$ is given by a smooth kernel $k(x,y)$ as a function of $x$ and $y$ off the diagonal, but which is typically singular on the diagonal. Such kernels arise, for instance, by applying a boundary integral formulation to second-order elliptic problems. Examples include the single-layer potential operator in Section 6.2 below or Abel-type operators with $k(x,y) = b(x,y)/|x-y|^\beta$ on $\Gamma = [0,1]$ for some $\beta > 0$ (see, e.g., Dahmen, Harbrecht and Schneider [6]). Assuming that $k$ is tractable only up to some experimental error, we postulate the knowledge of $dk_\delta(x,y) = dk(x,y) + \delta d\tilde{B}(x,y)$, where $\tilde{B}$ is a Brownian sheet. Assuming moreover that our data $g$ is perturbed by measurement noise as in (1.1), we recover our abstract framework.

*Statistical inference.* The widespread econometric model of instrumental variables (e.g., Hall and Horowitz [13]) is given by i.i.d. observations $(X_i, Y_i, W_i)$ for $i = 1, \ldots, n$, where $(X_i, Y_i)$ follow a regression model

$$Y_i = g(X_i) + U_i$$

with the exception that $\mathbb{E}[U_i|X_i] \neq 0$, but under the additional information given by the instrumental variables $W_i$ that satisfy $\mathbb{E}[U_i|W_i] = 0$. Denoting by $f_{XW}$ the joint density of $X$ and $W$, we define

$$k(x,z) := \int f_{XW}(x,w) f_{XW}(z,w) \, dw,$$

$$Kh(x) := \int k(x,z) h(z) \, dz.$$

To draw inference on $g$, we use the identity $Kg(x) = \mathbb{E}[\mathbb{E}[Y|W] f_{XW}(x,W)]$. The data easily allow estimation of the right-hand side and of the kernel function $k$. We face exactly an ill-posed inverse problem with errors in the operator, except for certain correlations between the two noise sources and for the fact that the noise is caused by a density estimation problem. Note that $K$ has a symmetric nonnegative kernel and is therefore self-adjoint and nonnegative on $L^2$. Hall and Horowitz [13] obtain in their Theorem 4.2 the linear rate of Section 3 when replacing their terms as follows: $\varepsilon = \delta = n^{-1/2}$, $t = \alpha$, $s = \beta + 1/2$, $d = 1$.

In other statistical problems random matrices or operators are of key importance or even the main subject of interest, for instance the linear response function in functional data analysis (e.g., Cai and Hall [2]) or the empirical covariance operator for stochastic processes (e.g., Reiss [17]).

*Numerical discretization.* Even if the operator is known, the numerical analyst is confronted with the same question of error in the operator under a different angle: up to which accuracy should the operator be discretized? Even more importantly, by not using all available information on the operator the objects typically have a sparse data structure and thus require far less memory and time of computation; see Dahmen, Harbrecht and Schneider [6].

**3. A linear estimation method.** In the following, we write $a \lesssim b$ when $a \leq cb$ for some constant $c > 0$ and $a \sim b$ when $a \lesssim b$ and $b \lesssim a$ simultaneously. The uniformity in $c$ will be obvious from the context.

3.1. *The linear Galerkin estimator.* We briefly study a linear projection estimator. Given $s > 0$ and $M > 0$, we first consider the Sobolev ball

$$W^s(M) := \{f \in H^s; \|f\|_{H^s} \leq M\}$$

as parameter space for the unknown $f$. Pick some $j \geq 0$ and let $V_j = \operatorname{span}\{\psi_\lambda, |\lambda| \leq j\}$ denote an approximation space associated with a $(\lfloor s \rfloor + 1)$-regular multiresolution analysis $(V_j)$; see Appendix A.6. We look for an estimator $\hat{f}_{\delta,\varepsilon} \in V_j$, solution to

(3.1) $$\langle K_\delta \hat{f}_{\delta,\varepsilon}, v \rangle = \langle g_\varepsilon, v \rangle \qquad \text{for all } v \in V_j.$$

This only makes sense if $K_\delta$ restricted to $V_j$ is invertible. We introduce the Galerkin projection (or *stiffness matrix*) of an operator $T$ onto $V_j$ by setting $T_j := P_j T|_{V_j}$, where $P_j$ is the orthogonal projection onto $V_j$, and set formally

(3.2) $$\hat{f}_{\delta,\varepsilon} := \begin{cases} K_{\delta,j}^{-1} P_j g_\varepsilon, & \text{if } \|K_{\delta,j}^{-1}\|_{V_j \to V_j} \leq \tau 2^{jt}, \\ 0, & \text{otherwise,} \end{cases}$$

where $\|T_j\|_{V_j \to V_j} = \sup_{v \in V_j, \|v\|_{L^2}=1} \|T_j v\|$ denotes the norm of the operator $T_j : (V_j, \|\bullet\|_{L^2}) \to (V_j, \|\bullet\|_{L^2})$. The estimator $\hat{f}_{\delta,\varepsilon}$ is specified by the level $j$ and the cut-off parameter $\tau > 0$ (and the choice of the multiresolution analysis).

3.2. *Result.* The ill-posedness comes from the fact that $K^{-1}$ is not $L^2$-continuous: we quantify the smoothing action by a *degree of ill-posedness* $t > 0$, which indicates that $K$ behaves roughly like $t$-fold integration. This is precisely defined by the following ellipticity condition in terms of the $L^2$-Sobolev norm $\|\bullet\|_{H^s}$ of regularity $s \in \mathbb{R}$; see Appendix A.6.

ASSUMPTION 3.1. $K$ is self-adjoint on $L^2(\mathcal{D})$, $K : L^2 \to H^t$ is continuous and $\langle K f, f \rangle \sim \|f\|_{H^{-t/2}}^2$.

As proved in Appendix A.6, Assumption 3.1 implies that the following "mapping constant" of $K$ with respect to the given multiresolution analysis $(V_j)$ is finite:

(3.3) $$Q(K) := \sup_{j \geq 0} 2^{-jt} \|K_j^{-1}\|_{V_j \to V_j} < \infty.$$

Introduce the integrated mean square error

$$\mathcal{R}(\hat{f}, f) := \mathbb{E}[\|\hat{f} - f\|_{L^2(\mathcal{D})}^2]$$

for an estimator $\hat{f}$ of $f$ and the rate exponent

$$r(s, t, d) := \frac{2s}{2s + 2t + d}.$$

PROPOSITION 3.2. *Let $Q > 0$. If the linear estimator $\hat{f}_{\delta,\varepsilon}$ is specified by $2^j \sim \max\{\delta, \varepsilon\}^{-2/(2s+2t+d)}$ and $\tau > Q$, then*

$$\sup_{f \in W^s(M)} \mathcal{R}(\hat{f}_{\delta,\varepsilon}, f) \lesssim \max\{\delta, \varepsilon\}^{2r(s,t,d)}$$

*holds uniformly over $K$ satisfying Assumption* 3.1 *with $Q(K) \leq Q$.*

The normalized rate $\max\{\delta, \varepsilon\}^{r(s,t,d)}$ gives the explicit interplay between $\varepsilon$ and $\delta$ and is indeed optimal over operators $K$ satisfying Assumption 3.1 with $Q(K) \leq Q$; see Section 5.2 below. Proposition 3.2 is essentially contained in Efromovich and Koltchinskii [11], but is proved in Section 7.1 as a central reference for the nonlinear results.

### 4. Two nonlinear estimation methods.

4.1. *Nonlinear Estimation* I. For $x > 0$ and two resolution levels $0 \leq j_0 < j_1$, define the level-dependent hard-thresholding operator $\mathcal{S}_x$ acting on $L^2(\mathcal{D})$ by

$$(4.1) \qquad \mathcal{S}_x(h) := \sum_{|\lambda| \leq j_1} \langle h, \psi_\lambda \rangle \mathbf{1}_{\{|\langle h, \psi_\lambda \rangle| \geq \kappa 2^{|\lambda| t} x \sqrt{(|\lambda| - j_0)_+}\}} \psi_\lambda,$$

for some constant $\kappa > 0$ and where $(\psi_\lambda)$ is a regular wavelet basis generating the multiresolution analysis $(V_j)$. Our first nonlinear estimator is defined by

$$(4.2) \qquad \hat{f}^I_{\delta,\varepsilon} := \mathcal{S}_{\max\{\delta,\varepsilon\}}(\hat{f}_{\delta,\varepsilon}),$$

where $\hat{f}_{\delta,\varepsilon}$ is the linear estimator (3.2) specified by the level $j_1$ and $\tau > 0$.

The factor $2^{|\lambda|t}$ in the threshold takes into account the increase in the noise level after applying the operator $K_{\delta,j_1}^{-1}$. The additional term $\sqrt{(|\lambda| - j_0)_+}$ is chosen to attain the exact minimax rate in the spirit of Delyon and Juditsky [8]. Hence, the nonlinear estimator $\hat{f}^I_{\delta,\varepsilon}$ is specified by $j_0$, $j_1$, $\tau$ and $\kappa$.

4.2. *Nonlinear Estimation* II. Our second method is conceptually different: we use matrix compression to remove the operator noise by thresholding $K_\delta$ in a first step and then apply the Galerkin inversion on the smoothed data $g_\varepsilon$. Let

$$(4.3) \qquad \hat{K}_\delta := \mathcal{S}^{\mathrm{op}}_\delta(K_{\delta,J}),$$

where $K_{\delta,J} = P_J K_\delta|_{V_J}$ is the Galerkin projection of the observed operator and $\mathcal{S}^{\mathrm{op}}_\delta$ is a hard-thresholding rule applied to the entries in the wavelet representation of the operator:

$$(4.4) \qquad T_J \mapsto \mathcal{S}^{\mathrm{op}}_\delta(T_J) := \sum_{|\lambda|,|\lambda'| \leq J} T_{\lambda,\lambda'} \mathbf{1}_{\{|T_{\lambda,\lambda'}| \geq \mathcal{T}(\delta)\}} \langle \bullet, \psi_\lambda \rangle \psi_{\lambda'},$$

where $\mathcal{T}(x) = \kappa x \sqrt{|\log x|}$ and $T_{\lambda,\lambda'} := \langle T\psi_\lambda, \psi_{\lambda'} \rangle$.

The estimator $\hat{g}_\varepsilon$ of the data is obtained by the classical hard-thresholding rule for noisy signals:

$$(4.5) \qquad \hat{g}_\varepsilon := \sum_{|\lambda| \leq J} \langle g_\varepsilon, \psi_\lambda \rangle \mathbf{1}_{\{|\langle g_\varepsilon, \psi_\lambda \rangle| \geq \mathcal{T}(\varepsilon)\}} \psi_\lambda.$$

After this preliminary step, we invert the linear system on the multiresolution space $V_J$ to obtain our second nonlinear estimator:

$$(4.6) \qquad \hat{f}_{\delta,\varepsilon}^{II} := \begin{cases} \hat{K}_\delta^{-1} \hat{g}_\varepsilon, & \|\hat{K}_\delta^{-1}\|_{H^t \to L^2} \leq \tau, \\ 0, & \text{otherwise.} \end{cases}$$

The nonlinear estimator $\hat{f}_{\delta,\varepsilon}^{II}$ is thus specified by $J$, $\kappa$ and $\tau$. Observe that this time we do not use level-dependent thresholds since we threshold the empirical coefficients directly.

## 5. Results for the nonlinear estimators.

### 5.1. *The setting.*

The nonlinearity of our two estimators permits to consider wider ranges of function classes for our target: we measure the smoothness $s$ of $f$ in $L^p$-norm, with $1 \leq p \leq 2$, in terms of Besov spaces $B_{p,p}^s$. The minimax rates of convergence are computed over Besov balls

$$V_p^s(M) := \{f \in B_{p,p}^s; \|f\|_{B_{p,p}^s} \leq M\}$$

with radius $M > 0$. We show that an elbow in the minimax rates is given by the critical line

$$(5.1) \qquad \frac{1}{p} = \frac{1}{2} + \frac{s}{2t+d},$$

considering $t$ and $d$ as fixed by the model setting. Equation (5.1) is linked to the geometry of inhomogeneous sparse signals that can be recovered in $L^2$-error after the action of $K$; see Donoho [9]. We retrieve the framework of Section 3 using $H^s = B_{2,2}^s$.

We prove in Section 5.2 that the first nonlinear estimator $\hat{f}_{\delta,\varepsilon}^{I}$ achieves the optimal rate over Besov balls $V_p^s(M)$. In Section 5.3 we further show that, under some mild restriction, the nonlinear estimator $\hat{f}_{\delta,\varepsilon}^{II}$ is adaptive in $s$ and nearly rate-optimal, losing a logarithmic factor in some cases.

### 5.2. *Minimax rates of convergence.*

In the following, we fix $s_+ \in \mathbb{N}$ and pick a wavelet basis $(\psi_\lambda)_\lambda$ associated with an $s_+$-regular multiresolution analysis $(V_j)$. The minimax rates of convergence are governed by the parameters $s \in (0, s_+)$, $p > 0$ and separate into two regions:

$$\text{dense region:} \qquad \mathcal{P}_{\text{dense}} := \left\{(s,p) : \frac{1}{p} < \frac{1}{2} + \frac{s}{2t+d}\right\},$$

$$\text{sparse region:} \qquad \mathcal{P}_{\text{sparse}} := \left\{(s,p) : \frac{1}{p} \geq \frac{1}{2} + \frac{s}{2t+d}\right\}.$$

It is implicitly understood that $B_{p,p}^s \subseteq L^2$ holds, that is, by Sobolev embeddings $s - d/p + d/2 \geq 0$. The terms dense and sparse refer to the form of the priors used to construct the lower bounds. Note that, an unavoidable logarithmic term appears in the sparse case.

THEOREM 5.1. *Let $Q > 0$. Specify the first nonlinear estimator $\hat{f}^I_{\delta,\varepsilon}$ by $2^{j_0} \sim \max\{\delta,\varepsilon\}^{-2/(2s+2t+d)}$, $2^{j_1} \sim \max\{\delta,\varepsilon\}^{-1/t}$, $\tau > Q$ and $\kappa > 0$ sufficiently large.*

- *For $(s,p) \in \mathcal{P}_{\mathrm{dense}}$ and $p \geq 1$ we have*

$$\sup_{f \in V^s_p(M)} \mathcal{R}(f^I_{\delta,\varepsilon}, f) \lesssim \max\{\delta,\varepsilon\}^{2r(s,t,d)},$$

*uniformly over $K$ satisfying Assumption 3.1 with $Q(K) \leq Q$.*
- *For $(s,p) \in \mathcal{P}_{\mathrm{sparse}}$ and $p \geq 1$ we have*

$$\sup_{f \in V^s_p(M)} \mathcal{R}(f^I_{\delta,\varepsilon}, f) \lesssim \max\{\delta\sqrt{|\log\delta|}, \varepsilon\sqrt{|\log\varepsilon|}\}^{2\tilde{r}(s,p,t,d)},$$

*uniformly over $K$ satisfying Assumption 3.1 with $Q(K) \leq Q$, where now*

$$\tilde{r}(s,p,t,d) := \frac{s + d/2 - d/p}{s + t + d/2 - d/p}.$$

A sufficient value for $\kappa$ can be made explicit by a careful study of Lemma 7.2 together with the proof of Delyon and Juditsky [8]; see the proofs below.

The rate obtained is indeed optimal in a minimax sense. The lower bound in the case $\delta = 0$ is classical (Nussbaum and Pereverzev [16]) and will not decrease for increasing noise levels $\delta$ or $\varepsilon$, whence it suffices to provide the case $\varepsilon = 0$.

The following lower bound can be derived from Efromovich and Koltchinskii [11] for $s > 0$, $p \in [1,\infty]$:

(5.2) $$\inf_{\hat{f}_\delta} \sup_{(f,K)\in\mathcal{F}_{s,p,t}} \mathcal{R}(\hat{f}_\delta, f) \gtrsim \delta^{2r(s,t,d)},$$

where the nonparametric class $\mathcal{F}_{s,p,t} = \mathcal{F}_{s,p,t}(M,Q)$ takes the form

$$\mathcal{F}_{s,p,t} = V^s_p(M) \times \{K \text{ satisfying Assumption 3.1 with } Q(K) \leq Q\}.$$

For $(s,p) \in \mathcal{P}_{\mathrm{dense}}$ the lower bound matches the upper bound attained by $\hat{f}^I_{\delta,\varepsilon}$. In Appendix A.5 we prove the following sparse lower bound:

THEOREM 5.2. *For $(s,p) \in \mathcal{P}_{\mathrm{sparse}}$ we have*

(5.3) $$\inf_{\hat{f}_\delta} \sup_{(K,f)\in\mathcal{F}_{s,p,t}} \mathcal{R}(\hat{f}_\delta, f) \gtrsim (\delta\sqrt{|\log\delta|})^{\tilde{r}(s,t,d)},$$

*and also the sparse rate of the first estimator $f^I_{\delta,\varepsilon}$ is optimal.*

5.3. *The adaptive properties of Nonlinear Estimation* II. We first state a general result which gives separate estimates for the two error levels of $\hat{f}^{II}_{\delta,\varepsilon}$ associated with $\delta$ and $\varepsilon$, respectively, leading to faster rates of convergence than in Theorem 5.1 in the case of sparse operator discretizations.

ASSUMPTION 5.3. $K: B^s_{p,p} \to B^{s+t}_{p,p}$ is continuous.

Furthermore, we state an ad hoc hypothesis on the sparsity of $K$. It is expressed in terms of the wavelet discretization of $K$ and is specified by parameters $(\bar{s}, \bar{p})$.

ASSUMPTION 5.4. For parameters $\bar{s} \geq 0$ and $\bar{p} > 0$ we have uniformly over all multi-indices $\lambda$

$$\|K\psi_\lambda\|_{B^{\bar{s}+t}_{\bar{p},\bar{p}}} \lesssim 2^{|\lambda|(\bar{s}+d/2-d/\bar{p})}.$$

Observe that this hypothesis follows from Assumption 5.3 with $(\bar{s}, \bar{p}) = (s, p)$, $\bar{p} \geq 1$, due to $\|\psi_\lambda\|_{B^s_{p,p}} \sim 2^{|\lambda|(s+d/2-d/p)}$. The case $\bar{p} < 1$, however, expresses high sparsity: if $K$ is diagonal in a regular wavelet basis with eigenvalues of order $2^{-|\lambda|t}$, then Assumption 5.4 holds for all $\bar{s}, \bar{p} \geq 0$. For a less trivial example of a sparse operator see Section 6.2. Technically, Assumption 5.4 will allow to control the error when thresholding the operator; see Proposition 7.4.

Finally, we need to specify a restriction on the linear approximation error expressed in terms of the regularity in $H^\alpha$:

$$(5.4) \quad \alpha \geq s\left(\frac{t+d}{s+t+d/2}\right)\min\left\{\frac{\log\varepsilon}{\log\delta}, 1\right\} \qquad \text{in the case } \delta > \varepsilon^{1+d/t}.$$

Then for $s \in (0, s_+)$, $p \geq 1$ and $\bar{p} > 0$ we obtain the following general result in the dense case.

THEOREM 5.5. *Grant Assumptions* 3.1, 5.3 *and* 5.4. *Let* $(s,p), (\bar{s}, \bar{p}) \in \mathcal{P}_{\text{dense}}$ *satisfy*

$$(5.5) \quad \frac{2\bar{s}+d-2d/\bar{p}}{2\bar{s}+2t+d} \leq \frac{2s-d}{2t+d} \qquad \text{with strict inequality for } p > 1$$

*and assume restriction (5.4) for $\alpha \geq 0$. Choose $\kappa > 0$ and $\tau > 0$ sufficiently large and specify $2^J \sim \min\{\varepsilon^{-1/t}, (\delta\sqrt{|\log\delta|})^{-1/(t+d)}\}$. Then*

$$\sup_{f \in V^s_p(M) \cap W^\alpha(M)} \mathcal{R}(\hat{f}^{II}_{\delta,\varepsilon}, f) \lesssim (\varepsilon\sqrt{|\log\varepsilon|})^{2r(s,t,d)} + (\delta\sqrt{|\log\delta|})^{2r(\bar{s},t,d)}.$$

The constant in the specification of $2^J$ cannot be too large; see the proof of Proposition 7.5. While the bounds for $\tau$ and $\kappa$ are explicitly computable from upper bounds on constants involved in the assumptions on the operator, they are in practice much too conservative, as is well known in the signal detection case (e.g., Donoho and Johnstone [10]) or the classical inverse problem case (Abramovich and Silverman [1]).

COROLLARY 5.6. *Grant Assumptions* 3.1 *and* 5.3. *Suppose* $(s,p) \in \mathcal{P}_{\mathrm{dense}}$ *and* $\alpha \geq 0$ *satisfies* (5.4). *Then*

$$\sup_{f \in V_p^s(M) \bigcap W^\alpha(M)} \mathcal{R}(\hat{f}_{\delta,\varepsilon}^{II}, f) \lesssim \max\{\varepsilon\sqrt{|\log \varepsilon|}, \delta\sqrt{|\log \delta|}\}^{2r(s,t,d)}$$

*follows from the smoothness restriction* $s \geq (d^2 + 8(2t+d)(d-d/p))^{1/2}/4$, *in particular in the cases* $p=1$ *or* $s > d(1 + \frac{2t}{d})^{1/2}/2$.

*If in addition* $d/p \leq d/2 + s(s-d/2)/(s+t+d/2)$ *holds, then we get rid of the linear restriction:*

$$\sup_{f \in V_p^s(M)} \mathcal{R}(\hat{f}_{\delta,\varepsilon}^{II}, f) \lesssim \max\{\varepsilon\sqrt{|\log \varepsilon|}, \delta\sqrt{|\log \delta|}\}^{2r(s,t,d)}.$$

PROOF. Set $\bar{s} = s$ and $\bar{p} = p$ and use that Assumption 5.3 implies Assumption 5.4. Then the smoothness restriction implies (5.5) and Theorem 5.5 applies. The particular cases follow because $s$ and $p$ are in $\mathcal{P}_{\mathrm{dense}}$.

By Sobolev embeddings, $B_{p,p}^s \subseteq W^\alpha$ holds for $s - d/p \geq \alpha - d/2$ and the last assertion follows by substituting in (5.4). □

We conclude that Nonlinear Estimation II attains the minimax rate up to a logarithmic factor in the dense case, provided the smoothness $s$ is not too small. For $(s,p) \in \mathcal{P}_{\mathrm{sparse}}$ the rate with exponent $\tilde{r}(s,p,t,d)$ is obtained via the Sobolev embedding $B_{p,p}^s \subseteq B_{\pi,\pi}^\sigma$ with $s - d/p = \sigma - d/\pi$ such that $(\sigma, \pi) \in \mathcal{P}_{\mathrm{dense}}$, and even exact rate-optimality follows in the sparse case.

## 6. Numerical implementation.

6.1. *Specification of the method.* While the mapping properties of the unknown operator $K$ along the scale of Sobolev or Besov spaces allow a proper mathematical theory and a general understanding, it is per se an asymptotic point of view: it is governed by the decay rate of the eigenvalues. For finite samples only the eigenvalues in the Galerkin projection $K_J$ matter, which will be close to the first $2^{Jd}$ eigenvalues of $K$. Consequently, even if the degree of ill-posedness of $K$ is known in advance (as is the case, e.g., in Reiss [17]), optimizing the numerical performance should rather rely on the induced norm $\|\bullet\|_{K_J} := \|K_J^{-1} \bullet\|_{L^2}$ on $V_J$ and not on $\|\bullet\|_{H^t}$.

Another practical point is that the cut-off rule using $\tau$ in the definitions (3.2) and (4.6) is not reasonable given just one sample, but needed to handle possibly highly distorted observations. An obvious way out is to consider only indices $J$ of the approximation space $V_J$ which are so small that $K_{\delta,J}$ remains invertible and not too ill-conditioned. Then the cut-off rule can be abandoned and the parameter $\tau$ is obsolete.

The estimator $\hat{f}_{\delta,\varepsilon}^{II}$ is therefore specified by choosing an approximation space $V_J$ and a thresholding constant $\kappa$. Since a thresholding rule is applied to both signal and operator, possibly different values of $\kappa$ can be used. In our experience, thresholds that are smaller than the theoretical bounds, but slightly larger than good choices in classical signal detection work well; see Abramovich and Silverman [1] for a similar observation.

The main constraint for selecting the subspace $V_J$ is that $J$ is not so large that $K_{\delta,J}$ is far away from $K_J$. By a glance at condition (7.6) in the proof of Theorem 5.5, working with $\|\bullet\|_{K_J}$ instead of $\|\bullet\|_{H^t}$ and with the observed operator before thresholding, we want that

$$\|K_{\delta,J} - K_J\|_{(V_J,\|\bullet\|_{L^2}) \to (V_J,\|\bullet\|_{K_J})} \leq \rho \|K_J^{-1}\|_{(V_J,\|\bullet\|_{K_J}) \to (V_J,\|\bullet\|_{L^2})}^{-1}$$

with some $\rho \in (0,1)$. This reduces to $\|(\mathrm{Id} - \delta K_{\delta,J}^{-1} \dot{B}_J)^{-1} - \mathrm{Id}\|_{V_J \to V_J} \leq \rho$, which by Lemma 7.1 is satisfied with very high probability provided

$$(6.1) \qquad \lambda_{\min}(K_{\delta,J}) \geq c\delta\sqrt{\dim(V_J)},$$

where $\lambda_{\min}(\bullet)$ denotes the minimal eigenvalue and $c > 0$ a constant depending on $\rho$ and the desired confidence. Based on these arguments we propose the following sequential data-driven rule to choose the parameter $J$:

$$(6.2) \qquad J := \min\{j \geq 0 | \lambda_{\min}(K_{\delta,j+1}) < c\delta \dim(V_{j+1})\}.$$

This rule might be slightly too conservative since after thresholding $\hat{K}_\delta$ will be closer to $K_J$ than $K_{\delta,J}$. It is, however, faster to implement and the desired confidence can be better tuned. In addition, a conservative choice of $J$ will only affect the estimation of very sparse and irregular functions.

6.2. *A numerical example.* We consider a single-layer logarithmic potential operator that relates the density of the electric charge on a cylinder of radius $r = 1/4$ to the induced potential on the same cylinder, when the cylinder is assumed to be infinitely long and homogeneous in that direction. Describing the angle by $e^{2\pi ix}$ with $x \in [0,1]$, the operator is given by

$$Kf(x) = \int_0^1 k(x,y)f(y)\,dy \qquad \text{with } k(x,y) = -\log(\tfrac{1}{2}|\sin(\pi(x-y))|).$$

The single-layer potential operator is known to satisfy a degree of ill-posedness $t = 1$ because of its logarithmic singularity on the diagonal. In Cohen, Hoffmann and Reiss [5] this operator has been used to demonstrate different solution methods for inverse problems with known operator: the singular-value decomposition (SVD), the linear Galerkin method and a nonlinear Galerkin method which corresponds to Nonlinear Estimation II in the case $\delta = 0$.

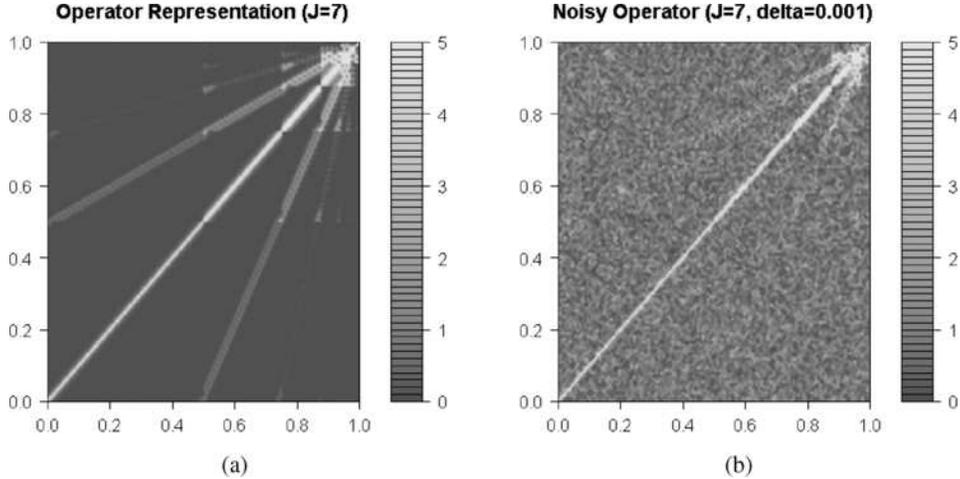

Fig. 1. *Wavelet representation of $K$* (a) *and $K_\delta$* (b).

The aim here is to compare the performance of the presented methods given that not $K$, but only a noisy version $K_\delta$ is available. Our focus is on the reconstruction properties under noise in the operator and we choose $\delta = 10^{-3}$, $\varepsilon = 10^{-5}$. As in Cohen, Hoffmann and Reiss [5] we consider the tent function

$$f(x) = \max\{1 - 30|x - \tfrac{1}{2}|, 0\}, \qquad x \in [0,1],$$

as object to be estimated. Its spike at $x = 1/2$ will be difficult to reconstruct.

For implementing the linear and the two nonlinear methods we use Daubechies wavelets of order 8 (with an extremal phase choice). We calculate the wavelet decomposition of $K$ and $f$ up to the scale $J_{\max} = 10$ by Mallat's pyramidal algorithm. For the nonlinear methods the large space $V_J$, on which the Galerkin inversion is performed, is determined by the rule (6.2) with $c = 5$. Figure 1(a) shows the modulus of the wavelet discretization $(|K_{\lambda,\mu}|)$ of the operator $K$ on $V_J$ with $J = 7$. Multi-indices with the same resolution level $j$ are presented next to each other; the resolution level $j$ decreases from left to right and from bottom to top. The units are multiples of $\delta$. The finger-like structure, showing large coefficients for low resolution levels, along the diagonal and certain subdiagonals, is typical for wavelet representations of integral (*Calderon–Zygmund*) operators and due to the support properties of the wavelets; see, for example, Dahmen, Harbrecht and Schneider [6]. In Figure 1(b) the modulus of the wavelet discretization of the noisy observation $K_\delta$ is shown. The structures off the main diagonal are hardly discernible.

The performance of the methods for this simulation setup are very stable for different noise realizations. In Figure 2(a) a typical linear estimation result for the choice $j = 5$ is shown along with the true function (dashed). Remark that because of $2^5/2^7 = 1/4$ the result is obtained by using only the values of $K_\delta$ that are depicted in the upper right quarter $[0.75, 1]^2$ of Figure 1(b). For the oracle choice $j = 5$ the root mean square error (RMSE) is minimal and evaluates to 0.029.

For the two nonlinear estimation methods, the approximation space $V_J$ (i.e., $V_{j_1}$ for Nonlinear Estimation I) chosen by the data-driven rule is $J = 7$ for all realizations. As to be expected, the simulation results deviate only marginally for different choices of $c \in [1, 20]$, giving either $J = 6$ or (mostly) $J = 7$. An implementation of Nonlinear Estimation I is based on a level-dependent thresholding factor which is derived from the average decay of the observed eigenvalues of $K_{\delta,J}$, ignoring the Delyon–Juditsky correction $\sqrt{(|\lambda| - j_0)_+}$. With the threshold base level $\kappa = 0.4$ (oracle choice) Nonlinear Estimation I produces an RMSE of 0.033. It shows a smaller error than the linear estimator at the flat parts far off the spike, but has difficulties with too large fluctuations close to the spike. The main underlying problem is that after the inversion the noise in the coefficients is heterogeneous even on the same resolution level which is not reflected by the thresholding rule.

Setting the base level $\kappa = 1.5$ for thresholding the operator and the data, the resulting estimator $f^{II}_{\delta,\varepsilon}$ of Nonlinear Estimation II is shown in Figure 2(b). It has by far the best performance among all three estimators with an RMSE of 0.022. The only artefacts, from an a posteriori perspective, are found next to the spike and stem from overlapping wavelets needed to reconstruct the spike itself.

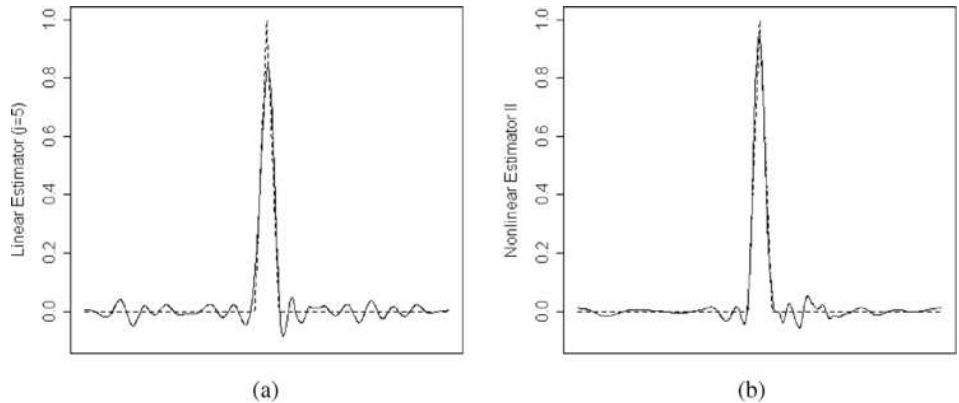

Fig. 2. *Linear estimator* (a) *and Nonlinear* II *estimator* (b).

In Cohen, Hoffmann and Reiss [5] simulations were performed for $\varepsilon = 2 \cdot 10^{-4}$ knowing the operator $K$ ($\delta = 0$). There the respective RMSE under oracle specifications is 0.024 (SVD), 0.023 (linear Galerkin), 0.019 (nonlinear Galerkin). In comparison we see that roughly the same accuracy is achieved in the case $\delta = 10^{-3}$, $\varepsilon = 10^{-5}$, which shows that the error in the operator is less severe than the error in the data. This observation is corroborated by further simulations for different values of $\delta$ and $\varepsilon$.

In order to understand why in this example the error in the operator is less severe and Nonlinear Estimation II performs particularly well, let us consider more generally the properties for thresholding a sparse operator representation as in Figure 1(a). This is exactly the point where Assumption 5.4 comes into play with $\bar{p} \in (0,1)$. To keep it simple, let us focus on the extreme case of an operator $K$ which is diagonalized by the chosen wavelet basis with eigenvalues $2^{-|\lambda|t}$. Then $K$ satisfies Assumption 5.4 for all $(\bar{s}, \bar{p})$ and by Theorem 5.5, choosing $\bar{p}$ such that $(\bar{s}, \bar{p}) \in \mathcal{P}_{\text{dense}}$ and restriction (5.5) is satisfied with equality, we infer

$$\sup_{f \in V_p^s(M) \cap W^\alpha(M)} \mathcal{R}(\hat{f}_{\delta,\varepsilon}^{II}, f) \lesssim (\varepsilon \sqrt{|\log \varepsilon|})^{2r(s,t,d)} + (\delta \sqrt{|\log \delta|})^{\min\{(2s-d)/t, 2\}}.$$

This rate is barely parametric in $\delta$ for not too small $s$. Hence, Nonlinear Estimation II can profit from the usually sparse wavelet representation of an operator, even without any specific tuning. This important feature is shared neither by Nonlinear Estimation I nor by the linear Galerkin method.

## 7. Proofs.

7.1. *Proof of Proposition* 3.2. By definition, $\mathcal{R}(\hat{f}_{\delta,\varepsilon}, f)$ is bounded by a constant times the sum of three terms $I + II + III$, where term $III$ comes from $\hat{f}_{\delta,\varepsilon} = 0$ if $\|K_{\delta,j}^{-1}\|_{V_j \to V_j} > \tau 2^{jt}$:

$$I := \|f - f_j\|_{L^2}^2,$$
$$II := \mathbb{E}[\|(K_{\delta,j}^{-1} P_j g_\varepsilon - f_j) \mathbf{1}_{\{\|K_{\delta,j}^{-1}\|_{V_j \to V_j} \leq \tau 2^{jt}\}}\|_{L^2}^2],$$
$$III := \|f\|_{L^2}^2 \mathbb{P}(\|K_{\delta,j}^{-1}\|_{V_j \to V_j} > \tau 2^{jt}).$$

*The term $I$.* This bias term satisfies under Assumption 3.1

$$\|f - f_j\|_{L^2}^2 \lesssim 2^{-2js} \sim \max\{\delta, \varepsilon\}^{4s/(2s+2t+d)}$$

by estimate (A.1) in the Appendix and thus has the right order.

*The term III.* For $\rho \in (0,1)$ let us introduce the event

(7.1) $$\Omega_{\rho,\delta,j} = \{\delta \|K_j^{-1}\dot{B}_j\|_{V_j \to V_j} \leq \rho\}.$$

On the event $\Omega_{\rho,\delta,j}$ the operator $K_{\delta,j} = K_j(\mathrm{Id} + \delta K_j^{-1}\dot{B}_j)$ is invertible with $\|K_{\delta,j}^{-1}\|_{V_j \to V_j} \leq (1-\rho)^{-1}\|K_j^{-1}\|_{V_j \to V_j}$ because

(7.2) $$\|(\mathrm{Id} + \delta K_j^{-1}\dot{B}_j)^{-1}\|_{V_j \to V_j} \leq \sum_{m \geq 0} \|\delta K_j^{-1}\dot{B}_j\|_{V_j \to V_j}^m \leq (1-\rho)^{-1}$$

follows from the usual Neumann series argument. By (3.3), the choice $\rho > 1 - Q/\tau \in (0,1)$ thus implies $\{\|K_{\delta,j}^{-1}\|_{V_j \to V_j} > \tau 2^{jt}\} \subseteq \Omega_{\rho,\delta,j}^c$. For $\eta = 1 - (2t+d)/(2s+2t+d) > 0$ and sufficiently small $\delta$, we claim that

(7.3) $$\mathbb{P}(\Omega_{\rho,\delta,j}^c) \leq \exp(-C\rho\delta^{-\eta}2^{2jd}) \qquad \text{for some } C > 0,$$

which implies that term *III* is of exponential order and hence negligible. To prove (7.3), we infer from (3.3)

$$\Omega_{\rho,\delta,j}^c \subseteq \{2^{-jd/2}\|\dot{B}_j\|_{V_j \to V_j} > \rho\delta^{-1}\|K_j^{-1}\|_{V_j \to V_j}^{-1}2^{-jd/2}\}$$
$$\subseteq \{2^{-jd/2}\|\dot{B}_j\|_{V_j \to V_j} > \rho\delta^{-1}Q^{-1}2^{-j(2t+d)/2}\}$$

and the claim (7.3) follows from $\delta^{-1}2^{-j(2t+d)/2} \gtrsim \delta^{-\eta}$ and the following classical bound for Gaussian random matrices:

LEMMA 7.1 ([7], Theorem II.4). *There are constants $\beta_0, c, C > 0$ such that*

$$\forall \beta \geq \beta_0 : \mathbb{P}(2^{-jd/2}\|\dot{B}_j\|_{V_j \to V_j} \geq \beta) \leq \exp(-c\beta^2 2^{2jd}),$$
$$\forall \beta \geq 0 : \mathbb{P}(2^{-jd/2}\|\dot{B}_j\|_{V_j \to V_j} \leq \beta) \leq (C\beta)^{2^{2jd}}.$$

*The term II.* Writing $P_j g_\varepsilon = P_j Kf + \varepsilon P_j \dot{W}$ and using the independence of the event $\Omega_{\rho,\delta,j}^c$ from $P_j \dot{W}$ (recall $\dot{B}$ and $\dot{W}$ are independent), we obtain

$$\mathbb{E}[\|K_{\delta,j}^{-1}P_j g_\varepsilon - f_j\|_{L^2}^2 \mathbf{1}_{\{\|K_{\delta,j}^{-1}\|_{V_j \to V_j} \leq \tau 2^{jt}\}} \mathbf{1}_{\Omega_{\rho,j,\delta}^c}]$$
$$\lesssim 2^{2jt}(\|P_j Kf\|_{L^2}^2 + \varepsilon^2 \mathbb{E}[\|P_j \dot{W}\|_{L^2}^2] + \|f_j\|_{L^2}^2)\mathbb{P}(\Omega_{\rho,\delta,j}^c).$$

Because of $\|P_j Kf\|_{L^2}^2 + \|f_j\|_{L^2}^2 \lesssim M^2$, $\mathbb{E}[\|P_j \dot{W}\|_{L^2}^2] \lesssim 2^{jd}$ and estimate (7.3), we infer that the above term is asymptotically negligible. Therefore, we are left with proving that $\mathbb{E}[\|K_{\delta,j}^{-1}P_j g_\varepsilon - f_j\|_{L^2}^2 \mathbf{1}_{\Omega_{\rho,j,\delta}}]$ has the right order. On $\Omega_{\rho,j,\delta}$ we consider the decomposition

(7.4) $$K_{\delta,j}^{-1}P_j g_\varepsilon - f_j = ((\mathrm{Id} + \delta K_j^{-1}\dot{B}_j)^{-1} - \mathrm{Id})f_j$$
$$+ \varepsilon(\mathrm{Id} + \delta K_j^{-1}\dot{B}_j)^{-1}K_j^{-1}P_j\dot{W}.$$

As for the second term on the right-hand side of (7.4), we have

$$\mathbb{E}[\varepsilon^2 \|(\mathrm{Id} + \delta K_j^{-1}\dot{B}_j)^{-1} K_j^{-1} P_j \dot{W}\|_{L^2}^2 \mathbf{1}_{\Omega_{\rho,\delta,j}}]$$
$$\leq \varepsilon^2 \mathbb{E}[\|(\mathrm{Id} + \delta K_j^{-1}\dot{B}_j)^{-1}\|_{V_j \to V_j}^2 \mathbf{1}_{\Omega_{\rho,\delta,j}}] \|K_j^{-1}\|_{V_j \to V_j}^2 \mathbb{E}[\|P_j \dot{W}\|_{L^2}^2]$$
$$\lesssim \varepsilon^2 2^{2jt} 2^{dj} \sim \max\{\delta, \varepsilon\}^{4s/(2s+2t+d)},$$

where we used again the independence of $\Omega_{\rho,\delta,j}$ and $P_j\dot{W}$ and the bound (7.2). The first term on the right-hand side of (7.4) is treated by

$$\mathbb{E}[\|\delta K_j^{-1}\dot{B}_j(\mathrm{Id} + \delta K_j^{-1}\dot{B}_j)^{-1} f_j\|_{L^2}^2 \mathbf{1}_{\Omega_{\rho,\delta,j}}]$$
$$\leq \delta^2 \|K_j^{-1}\|_{V_j \to V_j}^2 \|f_j\|_{L^2}^2 \mathbb{E}[\|\dot{B}_j\|_{V_j \to V_j}^2 \|(\mathrm{Id} + \delta K_j^{-1}\dot{B}_j)^{-1}\|_{V_j \to V_j}^2 \mathbf{1}_{\Omega_{\rho,\delta,j}}]$$
$$\lesssim \delta^2 \|K_j^{-1}\|_{V_j \to V_j}^2 \mathbb{E}[\|\dot{B}_j\|_{V_j \to V_j}^2]$$
$$\lesssim \delta^2 2^{2jt} 2^{dj} \lesssim \max\{\delta, \varepsilon\}^{4s/(2s+2t+d)},$$

where we successively used the triangle inequality, $\lim_{j\to\infty} \|f_j\|_{L^2} = \|f\|_{L^2}$ from (A.1), bound (7.2), Lemma 7.1 and (3.3).

### 7.2. Proof of Theorem 5.1.

*The main decomposition.* The error $\mathcal{R}(\hat{f}_{\delta,\varepsilon}^I, f)$ is bounded by a constant times the sum $I + II + III$ with

$$I := \|f - f_{j_1}\|_{L^2}^2,$$
$$II := \mathbb{E}[\|\mathcal{S}_{\max\{\delta,\varepsilon\}}(\hat{f}_{\delta,\varepsilon}) - f_{j_1}\|_{L^2}^2 \mathbf{1}_{\{\|K_{\delta,j_1}^{-1}\|_{V_{j_1} \to V_{j_1}} \leq \tau 2^{j_1 t}\}}],$$
$$III := \|f\|_{L^2}^2 \mathbb{P}(\|K_{\delta,j_1}^{-1}\|_{V_{j_1} \to V_{j_1}} > \tau 2^{j_1 t}).$$

For the term $I$, we use the bias estimate (A.1) and the choice of $2^{j_1}$. The term $III$ is analogous to the term $III$ in the proof of Proposition 7.1; we omit the details. To treat the main term $II$, we establish sharp error bounds for the empirical wavelet coefficients $\langle \hat{f}_{\delta,\varepsilon}, \psi_\lambda \rangle$ for $|\lambda| \leq j_1$.

*The empirical wavelet coefficients.* We consider again the event $\Omega_{\rho,\delta,j_1}$ from (7.1) with $j_1$ in place of $j$. On that event, we have the decomposition

$$\hat{f}_{\delta,\varepsilon} = K_{\delta,j_1}^{-1} P_{j_1} g_\varepsilon = f_{j_1} - \delta K_{j_1}^{-1}\dot{B}_{j_1} f_{j_1} + \varepsilon K_{j_1}^{-1} P_{j_1}\dot{W} + r_{\delta,j_1}^{(1)} + r_{\delta,j_1}^{(2)},$$

with

$$r_{\delta,j_1}^{(1)} = \sum_{n \geq 2}(-\delta K_{j_1}^{-1}\dot{B}_{j_1})^n f_{j_1},$$

$$r_{\delta,j_1}^{(2)} = -\varepsilon \delta K_{j_1}^{-1}\dot{B}_{j_1}(\mathrm{Id} + \delta K_{j_1}\dot{B}_{j_1})^{-1} K_{j_1}^{-1} P_{j_1}\dot{W}.$$

In the Appendix we derive the following properties.

LEMMA 7.2. *Let $|\lambda| \leq j_1$ and $\rho \in (0, 1 - Q/\tau)$. Under Assumption 3.1 the following decomposition holds:*

$$\delta \langle K_{j_1}^{-1} \dot{B}_{j_1} f_{j_1}, \psi_\lambda \rangle = \delta 2^{|\lambda| t} \|f_{j_1}\|_{L^2} c_\lambda \xi_\lambda,$$

$$\varepsilon \langle K_{j_1}^{-1} P_{j_1} \dot{W}, \psi_\lambda \rangle = \varepsilon 2^{|\lambda| t} \tilde{c}_\lambda \tilde{\xi}_\lambda,$$

$$\langle r_{\delta, j_1}^{(1)}, \psi_\lambda \rangle = \delta^2 2^{|\lambda| t} \|f_{j_1}\|_{L^2} 2^{j_1(t+d)} \zeta_{\lambda, j_1},$$

$$\langle r_{\delta, \varepsilon, j_1}^{(2)}, \psi_\lambda \rangle = \delta \varepsilon 2^{|\lambda| t} 2^{j_1(t+d/2)} \tilde{\zeta}_{\lambda, j_1},$$

*on $\Omega_{\rho, \delta, j_1}$, where $|c_\lambda|, |\tilde{c}_\lambda| \lesssim 1$, $\xi_\lambda$ and $\tilde{\xi}_\lambda$ are standard Gaussian variables and $\zeta_{\lambda, j_1}$, $\tilde{\zeta}_{\lambda, j_1}$ are random variables satisfying*

$$\max\{\mathbb{P}(\{|\zeta_{\lambda, j_1}| \geq \beta\} \cap \Omega_{\rho, \delta, j_1}), \mathbb{P}(\{|\tilde{\zeta}_{\lambda, j_1}| \geq \beta\} \cap \Omega_{\rho, \delta, j_1})\} \leq \exp(-c\beta 2^{2j_1 d})$$

*for all $\beta \geq \beta_0$ with some (explicitly computable) constants $\beta_0, c > 0$.*

From this explicit decomposition we shall derive the fundamental deviation bound

$$\begin{aligned}
(7.5) \quad &\mathbb{P}(\{2^{-|\lambda| t} |\langle \hat{f}_{\delta, \varepsilon}, \psi_\lambda \rangle - \langle f_{j_1}, \psi_\lambda \rangle| \geq \beta \max\{\delta, \varepsilon\}\} \cap \Omega_{\rho, \delta, j_1}) \\
&\leq 4 \exp(-C\beta \min\{\beta, 2^{j_1 d}\})
\end{aligned}$$

for all $|\lambda| \leq j_1$ and some explicitly computable constant $C > 0$. Once this is achieved, we are in the standard signal detection setting with exponentially tight noise. The asserted bound for term $II$ is then proved exactly following the lines in [8]; see also the heteroskedastic treatment in [14]. The only fine point is that we estimate the Galerkin projection $f_{j_1}$, not $f$, but by estimate (A.2) in the Appendix $\|f_{j_1}\|_{B_{p,p}^s} \lesssim \|f\|_{B_{p,p}^s}$.

It remains to establish the deviation bound (7.5). By Lemma 7.2, that probability is bounded by the sum of the four terms

$$P_I := \mathbb{P}\left(\|f_{j_1}\|_{L^2} c_\lambda |\xi_\lambda| \geq \frac{\beta}{4}\right),$$

$$P_{II} := \mathbb{P}\left(|\tilde{c}_\lambda \tilde{\xi}_\lambda| \geq \frac{\beta}{4}\right),$$

$$P_{III} := \mathbb{P}\left(\left\{\delta 2^{j_1(t+d)} \|f_{j_1}\|_{L^2} \zeta_{\lambda, j_1} \geq \frac{\beta}{4}\right\} \cap \Omega_{\rho, \delta, j_1}\right),$$

$$P_{IV} := \mathbb{P}\left(\left\{\delta 2^{j_1(t+d/2)} \tilde{\zeta}_{\lambda, j_1} \geq \frac{\beta}{4}\right\} \cap \Omega_{\rho, \delta, j_1}\right).$$

We obtain the bounds $P_I \leq \exp(-c_I \beta^2)$, $P_{II} \leq \exp(-c_{II} \beta^2)$ with some constants $c_I, c_{II} > 0$ by Gaussian deviations. The large deviation bound on $\zeta_{\lambda, j_1}$ in Lemma 7.2 implies with a constant $c_{III} > 0$

$$P_{III} \leq \exp(-c_{III} \beta 2^{-j_1(t+d-2d)} \delta^{-1}).$$

Equally, $P_{IV} \leq \exp(-c_{IV}\beta 2^{-j_1(t+d/2-2d)}\delta^{-1})$ follows, which proves (7.5) with some $C > 0$ depending on $c_I$ to $c_{IV}$ since $\delta^{-1} \gtrsim 2^{j_1 t}$ by construction.

7.3. *Proof of Theorem* 5.5. The proof of Theorem 5.5 is a combination of a deviation bound for the hard-thresholding estimator in $H^t$-loss together with an error estimate in operator norm. The following three estimates are the core of the proof and seem to be new. They are proved in the Appendix.

PROPOSITION 7.3 (Deviation in $H^t$-norm). *Assume* $\kappa > 4\sqrt{t/d}$, $2^J \lesssim \varepsilon^{-1/t}$ *and that* $s \geq 0$, $p > 0$ *are in* $\mathcal{P}_{\text{dense}}$. *Then there exist constants* $c_0, \eta_0, R_0 > 0$ *such that for all functions* $g \in B_{p,p}^{s+t}$ *the hard-thresholding estimate* $\hat{g}_\varepsilon$ *in* (4.5) *satisfies with* $m := \max\{\|P_J g\|_{B_{p,p}^{s+t}}, \|P_J g\|_{B_{p,p}^{s+t}}^{p/2}\}$:

$$\forall \eta \geq \eta_0 : \mathbb{P}(\mathcal{T}(\varepsilon)^{-r(s,t,d)}\|\hat{g}_\varepsilon - P_J g\|_{H^t} \geq \eta m^{1-r(s,t,d)}) \lesssim \varepsilon^{c_0 \eta^2} + \varepsilon^{\kappa^2/8 - d/t},$$

$$\forall R \geq R_0 : \mathbb{P}(\|\hat{g}_\varepsilon - P_J g\|_{H^t} \geq m + R) \lesssim \varepsilon^{\kappa^2/16 - d/t} R^{-4}.$$

PROPOSITION 7.4 (Estimation in operator norm, $L^2$-bound). *Suppose* $\kappa^2 \geq 32 \max\{d/t, 1 + d(2t+d)/(4t(t+d))\}$. *Grant Assumption* 5.4 *with* $\bar{s} > 0, \bar{p} > 0$ *satisfying restriction* (5.5). *Then*

$$\mathbb{E}[\|\hat{K}_\delta - K_J\|_{(V_J, \|\bullet\|_{B_{p,p}^s}) \to H^t}^2] \lesssim (\delta\sqrt{|\log \delta|})^{2r(\bar{s}, t, d)}.$$

PROPOSITION 7.5 (Estimation in operator norm, deviation bound). *Suppose* $K_\infty := \sup_{\mu, \lambda} 2^{|\lambda|t} |\langle K\psi_\mu, \psi_\lambda\rangle| < \infty$. *Then for all* $\eta > 0$

$$\mathbb{P}(\|\hat{K}_\delta - K_J\|_{(V_J, \|\bullet\|_{L^2}) \to H^t} \geq c_0 |\log \delta|^{-1/2} + \eta) \lesssim \delta^{\eta \min\{\kappa^2/2 - 2d/(t+d), 1/2q_1\}},$$

*with* $q_1 := 2^J(\delta\sqrt{|\log \delta|})^{1/(t+d)}$ *and a constant* $c_0$ *depending only on* $K_\infty$.

For $\rho \in (0,1)$ we introduce the event

$$(7.6) \quad \Omega_{\rho,\delta,J}^{II} := \{\|\hat{K}_\delta - K_J\|_{(V_J, \|\bullet\|_{L^2}) \to H^t} \leq \rho \|K_J^{-1}\|_{(V_J, \|\bullet\|_{H^t}) \to L^2}^{-1}\}.$$

The Neumann series representation implies that on $\Omega_{\rho,\delta,J}^{II}$ the random operator

$$\hat{K}_\delta : (V_J, \|\bullet\|_{L^2}) \to (V_J, \|\bullet\|_{H^t})$$

is invertible with norm $\|\hat{K}_\delta^{-1}\| \leq (1-\rho)^{-1} \|K_J^{-1}\|$. For the subsequent choice $\rho \in (0, 1 - \|K_J^{-1}\|/\tau)$ this bound is smaller than the cut-off value $\tau$. On $\Omega_{\rho,\delta,J}^{II}$ we bound $\|\hat{f}_{\delta,\varepsilon}^{II} - f\|_{L^2}$ by

$$\|\hat{K}_\delta^{-1}(\hat{g}_\varepsilon - P_J g)\|_{L^2} + \|(\hat{K}_\delta^{-1} - K_J^{-1})P_J g\|_{L^2} + \|f_J - f\|_{L^2}.$$

The last term is the bias and has the right order by estimate (A.1) and the restriction (5.4) on $\alpha$. The first two stochastic errors are further bounded by

$$\|\hat{K}_\delta^{-1}\|_{(V_J, \|\bullet\|_{H^t}) \to L^2}(\|\hat{g}_\varepsilon - P_J g\|_{H^t} + \|\hat{K}_\delta - K_J\|_{(V_J, \|\bullet\|_{B_{p,p}^s}) \to H^t}\|f_J\|_{B_{p,p}^s}).$$

Because of $\|\hat{K}_\delta^{-1}\|_{(V_J, \|\bullet\|_{H^t}) \to L^2} \leq \tau$, the assertion on $\Omega_{\rho,\delta,J}^{II}$ follows from the standard risk estimate in $H^t$-loss (cf. Proposition 7.3 or [14], Theorem 3.1)

$$\mathbb{E}[\|\hat{g}_\varepsilon - P_J g\|_{H^t}^2] \lesssim \mathcal{T}(\varepsilon)^{2r(s,t,d)},$$

from the operator norm estimate of Proposition 7.4 and $\|f_J\|_{B_{p,p}^s} \lesssim \|f\|_{B_{p,p}^s}$; see (A.2).

On the complement $(\Omega_{\rho,\delta,J}^{II})^c$ the risk of $\hat{f}_{\delta,\varepsilon}^{II}$, conditional on $\dot{B}$, is uniformly bounded thanks to the cut-off rule in the construction. Assumption 5.4 and the symmetry of $K$ imply $2^{|\lambda|t}|\langle K\psi_\mu, \psi_\lambda\rangle| \lesssim 2^{-||\mu|-|\lambda||(s+d/2-d/p)}$. Consequently, Proposition 7.5 is applicable and a sufficiently large choice of $\kappa$ and a sufficiently small choice of $q_1$ by means of an appropriate choice of the constant in the specification of $2^J$ give $\mathbb{P}((\Omega_{\rho,\delta,J}^{II})^c) \lesssim \delta^2$, which ends the proof.

## APPENDIX

### A.1. Proof of Lemma 7.2.

*First equality.* By Assumption 3.1, $K_{j_1}$ is symmetric and thus $\delta\langle K_{j_1}^{-1}\dot{B}_{j_1}f_{j_1}, \psi_\lambda\rangle = \delta\langle \dot{B}_{j_1}f_{j_1}, K_{j_1}^{-1}\psi_\lambda\rangle$. This is a centered Gaussian random variable with variance $\delta^2\|f_{j_1}\|_{L^2}^2\|K_{j_1}^{-1}\psi_\lambda\|_{L^2}^2$. Assumption 3.1 gives $\|K_{j_1}^{-1}\psi_\lambda\|_{L^2}^2 \lesssim \|\psi_\lambda\|_{H^t}^2 \lesssim 2^{2|\lambda|t}$ (see Appendix A.6) and the first equality follows from estimate (A.1).

*Second equality.* We write $\varepsilon\langle K_{j_1}^{-1}P_{j_1}\dot{W}, \psi_\lambda\rangle = \varepsilon\langle \dot{W}, K_{j_1}^{-1}\psi_\lambda\rangle$, which is centered Gaussian with variance $\varepsilon^2\|K_{j_1}^{-1}\psi_\lambda\|_{L^2}^2$, and the foregoing arguments apply.

*Third equality.* On $\Omega_{\rho,\delta,j_1}$ the term $|\langle r_{\delta,j_1}^{(1)}, \psi_\lambda\rangle|$ equals

$$|\langle(\delta K_{j_1}^{-1}\dot{B}_{j_1})^2(\mathrm{Id}+\delta K_{j_1}^{-1}\dot{B}_{j_1})^{-1}f_{j_1}, \psi_\lambda\rangle|$$
$$= \delta^2|\langle\dot{B}_{j_1}K_{j_1}^{-1}\dot{B}_{j_1}(\mathrm{Id}+\delta K_{j_1}^{-1}\dot{B}_{j_1})^{-1}f_{j_1}, K_{j_1}^{-1}\psi_\lambda\rangle|$$
$$\leq \delta^2\|\dot{B}_{j_1}\|_{V_{j_1}\to V_{j_1}}^2\|K_{j_1}^{-1}\|_{V_{j_1}\to V_{j_1}}\|(\mathrm{Id}+\delta K_{j_1}^{-1}\dot{B}_{j_1})^{-1}\|_{V_{j_1}\to V_{j_1}}$$
$$\quad \times \|f_{j_1}\|_{L^2}\|K_{j_1}^{-1}\psi_\lambda\|_{L^2}$$
$$\lesssim \delta^2\|\dot{B}_{j_1}\|_{V_{j_1}\to V_{j_1}}^2 2^{j_1 t}2^{|\lambda|t},$$

where we successively applied the Cauchy–Schwarz inequality, (3.3), (7.2) on $\Omega_{\rho,\delta,j_1}$ and (A.1) together with the same arguments as before to bound $\|K_{j_1}^{-1}\psi_\lambda\|_{L^2}$. Lemma 7.1 yields the result.

*Fourth equality.* Since $\dot{W}$ and $\dot{B}$ are independent, we have that, conditional on $\dot{B}$, the random variable $\langle r_{\delta,\varepsilon,j_1}^{(2)}, \psi_\lambda\rangle \mathbf{1}_{\Omega_{\rho,\delta,j_1}}$ is centered Gaussian with conditional variance

$$\delta^2\varepsilon^2\|(K_{j_1}^{-1}\dot{B}_{j_1}(\mathrm{Id}+\delta K_{j_1}^{-1}\dot{B}_{j_1})^{-1}K_{j_1}^{-1})^*\psi_\lambda\|_{L^2}^2\mathbf{1}_{\Omega_{\rho,\delta,j_1}}$$
$$= \delta^2\varepsilon^2\|(\dot{B}_{j_1}(\mathrm{Id}+\delta K_{j_1}^{-1}\dot{B}_{j_1})^{-1}K_{j_1}^{-1})^*K_{j_1}^{-1}\psi_\lambda\|_{L^2}^2\mathbf{1}_{\Omega_{\rho,\delta,j_1}}$$
$$\lesssim \delta^2\varepsilon^2\|(\dot{B}_{j_1}(\mathrm{Id}+\delta K_{j_1}^{-1}\dot{B}_{j_1})^{-1}K_{j_1}^{-1})^*\|_{V_{j_1}\to V_{j_1}}^2 2^{2|\lambda|t}\mathbf{1}_{\Omega_{\rho,\delta,j_1}}$$
$$\lesssim \delta^2\varepsilon^2 2^{2(|\lambda|+j_1)t}\|\dot{B}_{j_1}^*\|_{V_{j_1}\to V_{j_1}}^2$$

by (3.3) and estimate (7.2), which is not affected when passing to the adjoint, up to an appropriate modification of $\Omega_{\rho,\delta,j_1}$ incorporating $B_{j_1}^*$. We conclude by applying Lemma 7.1 which is also not affected when passing to the adjoint.

**A.2. Proof of Proposition 7.3.** Denote by $g^\lambda$ and $g_\varepsilon^\lambda$ the wavelet coefficients of $g$ and $g_\varepsilon$. We have

$$\|\hat{g}_\varepsilon - P_J g\|_{H^t}^2 \sim \sum_{|\lambda|\leq J} 2^{2|\lambda|t}(g_\varepsilon^\lambda \mathbf{1}_{\{|g_\varepsilon^\lambda|\geq \mathcal{T}(\varepsilon)\}} - g^\lambda)^2.$$

The usual decomposition yields a bound of the right-hand side by the sum of four terms $I + II + III + IV$ with

$$I := \sum 2^{2|\lambda|t}(g_\varepsilon^\lambda - g^\lambda)^2 \mathbf{1}_{\{|g^\lambda|\geq (1/2)\mathcal{T}(\varepsilon)\}},$$
$$II := \sum 2^{2|\lambda|t}(g_\varepsilon^\lambda - g^\lambda)^2 \mathbf{1}_{\{|g_\varepsilon^\lambda - g^\lambda|> (1/2)\mathcal{T}(\varepsilon)\}},$$
$$III := \sum 2^{2|\lambda|t}(g^\lambda)^2 \mathbf{1}_{\{|g_\varepsilon^\lambda - g^\lambda|> \mathcal{T}(\varepsilon)\}},$$
$$IV := \sum 2^{2|\lambda|t}(g^\lambda)^2 \mathbf{1}_{\{|g^\lambda|< 2\mathcal{T}(\varepsilon)\}},$$

and where the sums in $\lambda$ range through the set $\{|\lambda|\leq J\}$.

*The term IV.* This approximation term is bounded by

$$\sum_{j\leq J} 2^{2jt} \sum_{|\lambda|=j} (2\mathcal{T}(\varepsilon))^{2-p}\min\{(g^\lambda)^p, (2\mathcal{T}(\varepsilon))^p\}$$
$$\lesssim \mathcal{T}(\varepsilon)^{2-p}\sum_{j\leq J} 2^{2jt}\min\{\|P_J g\|_{B_{p,p}^{s+t}}^p 2^{-j(s+t+d/2-d/p)p}, 2^{jd}\mathcal{T}(\varepsilon)^p\}$$

which is of order $\mathcal{T}(\varepsilon)^2 2^{\bar{j}(2t+d)}$ with
$$2^{\bar{j}(2s+2t+d)} \sim \min\{\|P_J g\|^2_{B^{s+t}_{p,p}}\mathcal{T}(\varepsilon)^{-2}, 2^{J(2s+2t+d)}\}.$$

Therefore, we obtain $IV \lesssim \|P_J g\|^{2-2r(s,t,d)}_{B^{s+t}_{p,p}}\mathcal{T}(\varepsilon)^{2r(s,t,d)}$.

*The term I.* For this second approximation term we need to introduce the random variables
$$\xi_j := \frac{\varepsilon^{-2}}{\#\{|\lambda|=j, |g^\lambda| \geq (1/2)\mathcal{T}(\varepsilon)\}} \sum_{|\lambda|=j} (g^\lambda_\varepsilon - g^\lambda)^2 \mathbf{1}_{\{|g^\lambda| \geq (1/2)\mathcal{T}(\varepsilon)\}}.$$

Using $\mathbf{1}_{\{|g^\lambda| \geq (1/2)\mathcal{T}(\varepsilon)\}} \leq |2g^\lambda/\mathcal{T}(\varepsilon)|^p$, we obtain for the least favorable case $1/p = 1/2 + s/(2t+d)$ that term $I$ is bounded by
$$\sum_{j \leq J} 2^{2jt}\varepsilon^2 \xi_j \sum_{|\lambda|=j} \mathbf{1}_{\{|g^\lambda| \geq (1/2)\mathcal{T}(\varepsilon)\}}$$
$$\lesssim \sum_{j \leq J} 2^{2jt}\varepsilon^2 \xi_j \min\left\{\mathcal{T}(\varepsilon)^{-p} \sum_{|\lambda|=j} |g^\lambda|^p, 2^{jd}\right\}$$
$$\lesssim \sum_{j \leq J} \varepsilon^2 \xi_j \min\{\mathcal{T}(\varepsilon)^{-p} 2^{-j(s+t+d/2-d/p)p+2jt}\|P_J g\|^p_{B^{s+t}_{p,p}}, 2^{j(2t+d)}\}.$$

Now observe that, as for term $IV$, the following inequality holds:
$$\sum_{j \leq J} \varepsilon^2 \min\{\mathcal{T}(\varepsilon)^{-p} 2^{-j(s+t+d/2-d/p)p+2jt}\|P_J g\|^p_{B^{s+t}_{p,p}}, 2^{j(2t+d)}\} \sim \varepsilon^2 2^{\bar{j}(2t+d)}$$
$$\text{with } 2^{\bar{j}(2s+2t+d)} \sim \min\{\|P_J g\|^p_{B^{s+t}_{p,p}}\mathcal{T}(\varepsilon)^{-2}, 2^{J(2s+2t+d)}\}.$$

By definition, each $\xi_j$ has a normalized (to expectation 1) $\chi^2$-distribution and so has any convex combination $\sum_j a_j \xi_j$. For the latter we infer $\mathbb{P}(\sum_j a_j \xi_j \geq \eta^2) \leq e^{-\eta^2/2}$, $\eta \geq 1$, by regarding the extremal case of only one degree of freedom. Consequently, we obtain $\mathbb{P}(c_1 \varepsilon^{-2} 2^{-\bar{j}(2t+d)} I \geq \eta^2) \leq e^{-\eta^2/2}$ with some constant $c_1 > 0$. Substituting for $\bar{j}$, we conclude with another constant $c_2 > 0$ that
$$\mathbb{P}(I \geq \tfrac{1}{2}\eta^2 \|P_J g\|^p_{B^s_{p,p}}(\mathcal{T}(\varepsilon)\|P_J g\|^{-p/2}_{B^{s+t}_{p,p}})^{2r(s,t,d)}) \leq \exp(-c_2 \eta^2 |\log \varepsilon|) = \varepsilon^{c_2 \eta^2}.$$

*The terms II and III.* For these deviation terms we obtain by independence and a Gaussian tail estimate
$$\mathbb{P}(\{II = 0\} \cap \{III = 0\}) \geq \mathbb{P}(|g^\varepsilon_\lambda - g_\lambda| \leq \tfrac{1}{2}\mathcal{T}(\varepsilon) \text{ for all } |\lambda| \leq J)$$
$$\geq (1 - \exp(-\kappa^2 |\log \varepsilon|/8))^{\#V_J}.$$
Using $\#V_J \sim 2^{Jd} \lesssim \varepsilon^{-d/t}$, we derive $\mathbb{P}(II + III > 0) \lesssim \varepsilon^{\kappa^2/8 - d/t}$.

*The first assertion.* We obtain for some $\eta_0 \geq 1$ and all $\eta \geq \eta_0$:

$$\mathbb{P}(\|\hat{g}_\varepsilon - P_J g\|_{H^t} \geq \eta m^{1-r(s,td)} \mathcal{T}(\varepsilon)^{r(s,t,d)})$$
$$\leq \mathbb{P}(I > \tfrac{1}{2}\eta^2 \|P_J g\|_{B^{s+t}_{p,p}}^{p(1-r(s,t,d))} \mathcal{T}(\varepsilon)^{2r(s,t,d)}) + \mathbb{P}(II + III > 0)$$
$$+ \mathbb{P}(IV > \tfrac{1}{2}\eta^2 \|P_J g\|_{B^{s+t}_{p,p}}^{2-2r(s,t,d)} \mathcal{T}(\varepsilon)^{2r(s,t,d)})$$
$$\lesssim \varepsilon^{c_2 \eta^2} + \varepsilon^{\kappa^2/8 - d/t} + 0.$$

*The second assertion.* We show that the deviation terms are well bounded in probability. While obviously $III \leq \|P_J g\|_{H^t}^2 \lesssim \|P_J g\|_{B^{s+t}_{p,p}}^2$ holds,

$$\mathbb{E}[II] \leq \sum_{|\lambda| \leq J} 2^{2|\lambda|t} \mathbb{E}[(g_\varepsilon^\lambda - g^\lambda)^4]^{1/2} \mathbb{P}(|g_\varepsilon^\lambda - g^\lambda| > \mathcal{T}(\varepsilon)/2)^{1/2}$$

is bounded in order by $2^{J(2t+d)} \varepsilon^2 \exp(\kappa^2 |\log \varepsilon|/8)^{1/2} \sim \varepsilon^{\kappa^2/16 - d/t}$ due to $2^J \lesssim \varepsilon^{-1/t}$. In the same way we find

$$\mathrm{Var}[II] \leq \sum_{|\lambda| \leq J} 2^{4|\lambda|t} \mathbb{E}[(g_\varepsilon^\lambda - g^\lambda)^8]^{1/2} \mathbb{P}(|g_\varepsilon^\lambda - g^\lambda| > \mathcal{T}(\varepsilon)/2)^{1/2} \lesssim \varepsilon^{\kappa^2/16 - d/t}.$$

By Chebyshev's inequality, we infer $\mathbb{P}(II \geq R^2) \lesssim \varepsilon^{\kappa^2/16 - d/t} R^{-4}$ for $R > 0$. Since the above estimates of the approximation terms yield superoptimal deviation bounds, the estimate follows for sufficiently large $R$.

**A.3. Proof of Proposition 7.4.** The wavelet characterization of Besov spaces (cf. Appendix A.6) together with Hölder's inequality for $p^{-1} + q^{-1} = 1$ yields

$$\|\hat{K}_\delta - K_J\|_{(V_J, \|\bullet\|_{B^s_{p,p}}) \to H^t}$$
$$\sim \sup_{\|(a_\mu)\|_{\ell^p} = 1} \left\| (\hat{K}_\delta - K_J) \left( \sum_{|\mu| \leq J} 2^{-|\mu|(s+d/2-d/p)} a_\mu \psi_\mu \right) \right\|_{H^t}$$
$$\leq \|(2^{-|\mu|(s+d/2-d/p)} \|(\hat{K}_\delta - K_J) \psi_\mu\|_{H^t})_{|\mu| \leq J}\|_{\ell^q}$$
$$\leq \|(2^{-|\mu|(s+d/2-d/p)} \|K_J \psi_\mu\|_{B^{\bar{s}+t}_{\bar{p},\bar{p}}}^{1-r(\bar{s},t,d)})_{|\mu| \leq J}\|_{\ell^q}$$
$$\times \sup_{|\mu| \leq J} \|K_J \psi_\mu\|_{B^{\bar{s}+t}_{\bar{p},\bar{p}}}^{r(\bar{s},t,d)-1} \|(\hat{K}_\delta - K_J) \psi_\mu\|_{H^t}.$$

Due to Assumption 5.4 the last $\ell^q$-norm can be estimated in order by

$$\|(2^{j(-(s-d/2) + (\bar{s}+d/2-d/\bar{p})(1-r(\bar{s},t,d)))})_{j \leq J}\|_{\ell^q},$$

which is of order 1 whenever restriction (5.5) is fulfilled.

By construction, $\hat{K}_\delta \psi_\mu$ is the hard-thresholding estimator for $K_J \psi_\mu$ given the observation of $K_{\delta,J} \psi_\mu$, which is $K_J \psi_\mu$ corrupted by white noise of level $\delta$. Therefore Proposition 7.3 applied to $K\psi_\mu$ and $\delta$ gives for any $\eta \geq \eta_0$:

$$\mathbb{P}(\|K_J\psi_\mu\|_{B_{\bar{p},\bar{p}}^{\bar{s}+t}}^{r(\bar{s},t,d)-1}\|(\hat{K}_\delta - K_J)\psi_\mu\|_{H^t} \geq \eta \mathcal{T}(\delta)^{r(\bar{s},t,d)}) \lesssim \delta^{c_0 \eta^2} + \delta^{\kappa^2/8 - d/t}.$$

By estimating the probability of the supremum by the sum over the probabilities, we obtain from above with a constant $c_1 > 0$ for all $\eta \geq \eta_0$:

$$\mathbb{P}(\|\hat{K}_\delta - K_J\|_{(V_J, \|\bullet\|_{B_{p,p}^s}) \to H^t} \geq \eta \mathcal{T}(\delta)^{r(\bar{s},t,d)})$$

$$\leq \sum_{|\mu| \leq J} \mathbb{P}(\|K\psi_\mu\|_{B_{\bar{p},\bar{p}}^{\bar{s}+t}}^{r(\bar{s},t,d)-1} \|(\hat{K}_\delta - K_J)\psi_\mu\|_{H^t} \geq c_1 \eta \mathcal{T}(\delta)^{r(\bar{s},t,d)})$$

$$\lesssim 2^{Jd}(\delta^{c_0 \eta^2} + \delta^{\kappa^2/8 - d/t})$$

$$\lesssim \delta^{c_0 \eta^2 - d/(t+d)} + \delta^{\kappa^2/8 - d(2t+d)/(t(t+d))}.$$

For a sufficiently large $\eta_1 > \eta_0$, depending only on $c_0$, $d$ and $t$, with $\gamma := \kappa^2/8 - d(2t+d)/(t(t+d)) > 0$, we thus obtain

$$\mathbb{P}(\|\hat{K}_\delta - K_J\|_{(V_J, \|\bullet\|_{B_{p,p}^s}) \to H^t} \geq \eta_1 \mathcal{T}(\delta)^{r(\bar{s},t,d)}) \lesssim \delta^\gamma.$$

By the above bound on the operator norm and Hölder's inequality for $q := \gamma/2 \geq 2$ and $\rho^{-1} + q^{-1} = 1$ together with the second estimate in Proposition 7.3, we find for some constant $R_0 > 0$:

$$\mathbb{E}[\|\hat{K}_\delta - K_J\|_{(V_J, \|\bullet\|_{B_{p,p}^s}) \to H^t}^2 \mathbf{1}_{\{\|\hat{K}_\delta - K_J\|_{(V_J, \|\bullet\|_{B_{p,p}^s}) \to H^t} \geq \eta_1 \mathcal{T}(\delta)^{r(\bar{s},t,d)}\}}]$$

$$\lesssim \mathbb{E}[\|\hat{K}_\delta - K_J\|_{(V_J, \|\bullet\|_{B_{p,p}^s}) \to H^t}^{2\rho}]^{1/\rho} \delta^{\gamma/q}$$

$$\leq \left(\int_0^\infty R^{2\rho-1} \mathbb{P}(\|\hat{K}_\delta - K_J\|_{(V_J, \|\bullet\|_{B_{p,p}^s}) \to H^t} \geq R) \, dR\right)^{1/\rho} \delta^2$$

$$\lesssim \left(R_0 + \int_{R_0}^\infty R^{2\rho-1} 2^{Jd} \delta^{\kappa^2/16 - d/t} R^{-4} \, dR\right)^{1/\rho} \delta^2$$

$$\lesssim \max\{\delta^{(\kappa^2/16 - 2d/t)/\rho}, 1\} \delta^2$$

which is of order $\delta^2$ by assumption on $\kappa$ and the assertion follows.

**A.4. Proof of Proposition 7.5.** For $|\mu|, |\lambda| \leq J$ we have for the entries in the wavelet representation

$$|(\hat{K}_\delta)_{\mu,\lambda} - K_{\mu,\lambda}| = |K_{\mu,\lambda}| \mathbf{1}_{\{|(K_\delta)_{\mu,\lambda}| \leq \mathcal{T}(\delta)\}} + \delta |\dot{B}_{\mu,\lambda}| \mathbf{1}_{\{|(K_\delta)_{\mu,\lambda}| > \mathcal{T}(\delta)\}}.$$

A simple rough estimate yields

$$|(\hat{K}_\delta)_{\mu,\lambda} - K_{\mu,\lambda}| \leq 2\mathcal{T}(\delta) + |K_{\mu,\lambda}| \mathbf{1}_{\{|(K_\delta - K)_{\mu,\lambda}| \geq \mathcal{T}(\delta)\}} + \delta |\dot{B}_{\mu,\lambda}|.$$

We bound the operator norm by the corresponding Hilbert–Schmidt norm and use $K_\infty < \infty$ to obtain

$$\|\hat{K}_\delta - K_J\|^2_{(V_J, \|\bullet\|_{L^2}) \to H^t}$$
$$\leq \sum_{|\mu|,|\lambda| \leq J} 2^{2|\lambda|t}((\hat{K}_\delta)_{\mu,\lambda} - K_{\mu,\lambda})^2$$
$$\lesssim 2^{2J(t+d)}\mathcal{T}(\delta)^2 + \#\{\delta | (K_\delta - K)_{\mu,\lambda}| \geq \mathcal{T}(\delta)\} + \delta^2 2^{2Jt} \sum_{|\mu|,|\lambda| \leq J} \dot{B}^2_{\mu,\lambda},$$

where the cardinality is taken for multi-indices $(\lambda, \mu)$ such that $|\lambda|, |\mu| \leq J$. The first term is of order $|\log \delta|^{-1}$. In view of $(K_\delta - K)_{\mu,\lambda} = \delta \dot{B}_{\mu,\lambda}$, the second term is a binomial random variable with expectation $2^{2Jd}\mathbb{P}(|\dot{B}_{\mu,\lambda}| \geq \kappa|\log\delta|^{1/2}) \lesssim \delta^{-2d/(t+d)+\kappa^2/2}$. An exponential moment bound for the binomial distribution yields

$$\mathbb{P}(\#\{|\dot{B}_{\mu,\lambda}| \geq \kappa\sqrt{|\log\delta|}\} \geq \eta) \lesssim \delta^{\eta(\kappa^2/2 - 2d/(t+d))}.$$

For the last term, we use an exponential bound for the deviations of a normalized $\chi^2$-distribution, as before, to infer from $2^{J(t+d)} \lesssim \mathcal{T}(\delta)$ that

$$\mathbb{P}\left(\delta^2 2^{2Jt} \sum_{|\mu|,|\lambda| \leq J} \dot{B}^2_{\mu,\lambda} \geq \eta\right) \leq \exp(-2^{-2J(t+d)-1}\delta^{-2}\eta) \leq \delta^{\eta/2q_1}$$

holds, which gives the result.

**A.5. Proof of Theorem 5.2.** To avoid singularity of the underlying probability measures we only consider the subclass $\mathcal{F}_0$ of parameters $(K, f)$ such that $Kf = y_0$ for some fixed $y_0 \in L^2$, that is, $\mathcal{F}_0 := \{(K,f) | f = K^{-1}y_0, K \in \mathcal{K}\}$, where $\mathcal{K} = \mathcal{K}_t(C)$ abbreviates the class of operators under consideration. We shall henceforth keep $y_0$ fixed and refer to the parameter $(K, f)$ equivalently just by $K$.

The likelihood $\Lambda(\bullet)$ of $\mathbb{P}^{K^2}$ under the law $\mathbb{P}^{K^1}$ corresponding to the parameters $K^i$, $i = 1, 2$, is

$$\Lambda(K^2, K^1) = \exp(\delta^{-1}\langle K^2 - K^1, \dot{B}\rangle_{\mathrm{HS}} - \tfrac{1}{2}\delta^{-2}\|K^1 - K^2\|^2_{\mathrm{HS}})$$

in terms of the scalar product and norm of the Hilbert space $HS(L^2)$ of Hilbert–Schmidt operators on $L^2$ and with a Gaussian white noise operator $\dot{B}$. In particular, the Kullback–Leibler divergence between the two measures equals $\tfrac{1}{2}\delta^{-2}\|K^1 - K^2\|^2_{\mathrm{HS}}$ and the two models remain contiguous for $\delta \to 0$ as long as the Hilbert–Schmidt norm of the difference remains of order $\delta$.

Let us fix the parameter $f_0 = \psi_{-1,0} = \mathbf{1}$ and the operator $K^0$ which, in a wavelet basis $(\psi_\lambda)_\lambda$, has diagonal form $K^0 = \mathrm{diag}(2^{-(|\lambda|+1)t})$. Then $K^0$ is ill-posed of degree $t$ and trivially obeys all the mapping properties imposed. Henceforth, $y_0 := K^0 f_0 = \mathbf{1}$ remains fixed.

For any $k = 0, \ldots, 2^{Jd} - 1$, introduce the symmetric perturbation $H^\varepsilon = (H^\varepsilon_{\lambda,\mu})$ with vanishing coefficients except for $H^\varepsilon_{(0,0),(J,k)} = 1$ and $H^\varepsilon_{(J,k),(0,0)} = 1$. Put $K^\varepsilon = K^0 + \gamma H^\varepsilon$ for some $\gamma > 0$. By setting $\gamma := \delta J$ we enforce $\|K^\varepsilon - K^0\|_{\mathrm{HS}} = \delta J$. For $f_\varepsilon := (K^\varepsilon)^{-1} y_0$, we obtain

$$f_\varepsilon - f_0 = ((K^\varepsilon)^{-1} - (K^0)^{-1}) y_0$$
$$= \gamma (K^\varepsilon)^{-1} H^\varepsilon f_0$$
$$= \gamma (K^\varepsilon)^{-1} \psi_{J,0}.$$

Now observe that $H^\varepsilon$ trivially satisfies the conditions

$$\tfrac{1}{2} |\langle H^\varepsilon f, f \rangle| \leq 2^{Jt} \|f\|^2_{H^{-t/2}}, \qquad \tfrac{1}{2} \|H^\varepsilon\|_{L^2 \to H^t} \leq 2^{Jt},$$
$$\tfrac{1}{2} \|H^\varepsilon\|_{B^s_{p,p} \to B^{s+t}_{p,p}} \leq 2^{J(t+s+d/2-d/p)}.$$

This implies that for $\gamma 2^{J(t+s+d/2-d/p)}$ sufficiently small $K^\varepsilon$ inherits the mapping properties from $K^0$. Hence,

$$\|f_\varepsilon - f_0\|_{L^2} \sim \gamma \|\psi_{J,0}\|_{H^t} = \gamma 2^{Jt},$$
$$\|f_\varepsilon - f_0\|_{B^s_{p,p}} \sim \gamma \|\psi_{J,0}\|_{B^{s+t}_{p,p}} = \gamma 2^{J(t+s+d/2-d/p)}$$

follows. In order to apply the classical lower bound proof in the sparse case ([15], Theorem 2.5.3) and thus to obtain the logarithmic correction, we nevertheless have to show that $f_\varepsilon - f_0$ is well localized. Using the fact that $((H^\varepsilon)^2)_{\lambda,\mu} = 1$ holds for coordinates $\lambda = \mu = (0,0)$ and $\lambda = \mu = (J,k)$, but vanishes elsewhere, we infer from the Neumann series representation

$$f_\varepsilon - f_0 = \sum_{m=1}^\infty (-\gamma H^\varepsilon)^m f_0$$
$$= \sum_{n=1}^\infty \gamma^{2n} f_0 - \sum_{n=0}^\infty \gamma^{2n+1} \psi_{J,k}$$
$$= \frac{\gamma}{1-\gamma^2} (\gamma f_0 - \psi_{J,k}).$$

Consequently, the asymptotics for $\gamma \to 0$ are governed by the term $-\gamma \psi_{J,k}$, which is well localized. The choice $2^J < \gamma^{-1/(t+s+d/2-d/p)}$ ensures that $\|f_\varepsilon\|_{B^s_{p,p}}$ remains bounded and we conclude by usual arguments; see Chapter 2 in [15] or the lower bound in [17].

**A.6. Some tools from approximation theory.** The material gathered here is classical; see, for example, [4]. We call a multiresolution analysis on $L^2(\mathcal{D})$ an increasing sequence of subspaces $(V_J)_{J\geq 0}$ generated by orthogonal wavelets $(\psi_\lambda)_{|\lambda|\leq J}$, where the multi-index $\lambda = (j,k)$ comprises the resolution level $|\lambda| := j \geq 0$ and the $d$-dimensional location parameter $k$. We use the fact that for regular domains $\#V_J \sim 2^{Jd}$ and denote the $L^2$-orthogonal projection onto $V_J$ by $P_J$.

Given an $s_+$-regular multiresolution analysis, $s_+ \in \mathbb{N}$, an equivalent norming of the Besov space $B_{p,p}^s$, $s \in (-s_+, s_+)$, $p > 0$, is given in terms of weighted wavelet coefficients:

$$\|f\|_{B_{p,p}^s} \sim \left( \sum_{j=-1}^\infty 2^{j(s+d/2-d/p)p} \sum_k |\langle f, \psi_{jk}\rangle|^p \right)^{1/p}.$$

For $p < 1$ the Besov spaces are only quasi-Banach spaces, but still coincide with the corresponding nonlinear approximation spaces; see Section 30 in [4]. If $s$ is not an integer or if $p=2$, the space $B_{p,p}^s$ equals the $L^p$-Sobolev space, which for $p=2$ is denoted by $H^s$. The Sobolev embedding generalizes to $B_{p,p}^s \subseteq B_{p',p'}^{s'}$ for $s \geq s'$ and $s - \frac{d}{p} \geq s' - \frac{d}{p'}$.

Direct and inverse estimates are the main tools in approximation theory. Using the equivalent norming, they are readily obtained for any $-s_+ < s' \leq s < s_+$:

$$\inf_{h_j \in V_j} \|f - h_j\|_{B_{p,p}^{s'}} \lesssim 2^{-(s-s')j} \|f\|_{B_{p,p}^s},$$

$$\forall h_j \in V_j : \|h_j\|_{B_{p,p}^s} \lesssim 2^{(s-s')j} \|h_j\|_{B_{p,p}^{s'}}.$$

In [5] it is shown that under Assumption 3.1 $\|K_j^{-1}\|_{H^t \to L^2} \lesssim 1$ and we infer from an inverse estimate $\|K_j^{-1}\|_{L^2 \to L^2} \lesssim 2^{jt}$.

Let us finally bound $\|f - f_j\|$ and $\|f_j\|$ for diverse norms. By definition, $f_j$ is the orthogonal projection of $f$ onto $V_j$ with respect to the scalar product $\langle K\cdot,\cdot\rangle$ such that $\|K^{1/2}(f-f_j)\|_{L^2} \leq \|K^{1/2}(\mathrm{Id}-P_j)f\|_{L^2}$ and by Assumption 3.1 $\|f - f_j\|_{H^{-t/2}} \lesssim \|(\mathrm{Id}-P_j)f\|_{H^{-t/2}}$. Using the equivalent (weighted) $\ell^p$-norms, we find $\|f - f_j\|_{B_{p,p}^{-t/2}} \lesssim \|(\mathrm{Id}-P_j)f\|_{B_{p,p}^{-t/2}}$ for any $p$. By a direct estimate, we obtain $\|(\mathrm{Id}-P_j)f\|_{B_{p,p}^{-t/2}} \lesssim 2^{-j(s+t/2)}\|f\|_{B_{p,p}^s}$ and

$$\|f - f_j\|_{B_{p,p}^{-t/2}} \lesssim 2^{-j(s+t/2)} \|f\|_{B_{p,p}^s},$$

hence

$$\|P_j f - f_j\|_{B_{p,p}^{-t/2}} \lesssim 2^{-j(s+t/2)} \|f\|_{B_{p,p}^s}.$$

An inverse estimate, applied to the latter inequality, yields together with the Sobolev embeddings ($p \leq 2$)

$$
\begin{aligned}
\|f - f_j\|_{L^2} &\leq \|f - P_j f\|_{L^2} + \|P_j f - f_j\|_{L^2} \\
&\lesssim 2^{-j(s+d/2-d/p)} \|f\|_{B^s_{p,p}}.
\end{aligned}
\tag{A.1}
$$

Merely an inverse estimate yields the stability estimate

$$
\|f_j\|_{B^s_{p,p}} \leq \|f_j - P_j f\|_{B^s_{p,p}} + \|P_j f\|_{B^s_{p,p}} \lesssim \|f\|_{B^s_{p,p}}.
\tag{A.2}
$$

**Acknowledgment.** Helpful comments by two anonymous referees are gratefully acknowledged.

CNRS-UMR 8050
AND
UNIVERSITY OF MARNE-LA-VALLÉE
FRANCE
E-MAIL: marc.hoffmann@univ-mlv.fr

INSTITUTE OF APPLIED MATHEMATICS
UNIVERSITY OF HEIDELBERG
69120 HEIDELBERG
GERMANY
E-MAIL: reiss@statlab.uni-heidelberg.de